# A Modification of Sufficient Conditions to Ensure the Exact Conic Relaxation

Tao Ding, *Member, IEEE*, Bo Zeng, *Member, IEEE*, Rui Bo, *Senior Member, IEEE*

*Abstract*—To solve the AC optimal power flow problem, it is proposed in [1, 2] that a convex conic approximation to branch flow model (BFM) can be obtained if we first eliminate phase angles of voltages and currents and then relax a set of equality constraints to second order conic ones. In particular, provided a set of sufficient conditions are satisfied, the conic relaxation is exact. We note, however, that those conditions do not always guarantee the exactness. In this letter, we analyze the argument of exactness and include a new condition that there is no line with negative reactance to ensure the conic formulation's exactness.

## I. INTRODUCTION

It is widely recognized that AC optimal power flow (ACOPF) problem is a very challenging non-convex optimization problem. One strategy that is computationally promising is to derive a strong convex relaxation that can be solved efficiently, and then to employ an effective method to convert a solution of that relaxation into a feasible one, which should be of a high quality to the original ACOPF problem. For example, a recent study in [1, 2] present a two-step relaxation method to obtain a second order cone program (SOCP) of the branch flow model (BFM) of ACOPF. Specifically, the first step is to eliminate restrictions of voltage and current angles; and the second step is to relax a set of quadratic equalities to their conic counterparts. As a result, these two steps lead to a SOCP problem that can be computed efficiently. In particular, authors of [1, 2] argue that, provided a set of sufficient conditions are satisfied, the second step is exact. Hence, they conclude that an optimal solution of the SOCP problem is also an optimal solution to the problem without restrictions from voltage and current angles. Therefore, as long as those angles can be recovered, that solution is optimal to the original ACOPF problem. For a radial network, angle recovery can be achieved easily. For a mesh network, the angle recovery condition might not be satisfied. However, by placing phase shifters on branches, the network can be convexified to satisfy that condition and then its OPF solution can be recovered.

In this letter, we concentrate on the conic relaxation, where the non-convex power flow equations (1) is relaxed to (2) by eliminating the voltage and current angles.

$$\begin{cases} P_i^g - P_i^c = U_i \sum U_j (G_{ij}\cos\theta_{ij} + B_{ij}\sin\theta_{ij}) & \forall i \in \{\mathbb{B}\} \\ Q_i^g - Q_i^c = U_i \sum U_j (G_{ij}\sin\theta_{ij} - B_{ij}\cos\theta_{ij}) & \forall i \in \{\mathbb{B}\} \end{cases} \quad (1)$$

where $U_i$ and $\theta_i$ are the voltage magnitude and angle of bus $i$; $P_i^g$ and $Q_i^g$ ($P_i^c$ and $Q_i^c$) are the active and reactive power genera-tion (load) at bus $i$; $G$ and $B$ are the real and imaginary parts of the admittance matrix; $\theta_{ij}=\theta_i-\theta_j$ is the angle difference between $i$ and $j$; $\{\mathbb{B}\}$ is the set of buses.

$$P_i^g - P_i^c = \sum_{k:i\to k} F_{ik} - \sum_{j:j\to i}(F_{ji} - r_{ji}l_{ji}), \forall i \in \{\mathbb{B}\} \quad (2\text{-a})$$

$$Q_i^g - Q_i^c = \sum_{k:i\to k} H_{ik} - \sum_{j:j\to i}(H_{ji} - x_{ji}l_{ji}), \forall i \in \{\mathbb{B}\} \quad (2\text{-b})$$

$$v_j = v_i - 2(r_{ij}F_{ij} + x_{ij}H_{ij}) + (r_{ij}^2 + x_{ij}^2)l_{ij}, \forall ij \in \{\mathbb{L}\} \quad (2\text{-c})$$

$$l_{ij}v_i = F_{ij}^2 + H_{ij}^2, \quad ij \in \{\mathbb{L}\} \quad (2\text{-d})$$

where $F_{ij}$ and $H_{ij}$ are the active and reactive branch flows from $i$ to $j$; $r_{ij}$ and $x_{ij}$ are the resistance and reactance of line $ij$; $v_i=U_i^2$; $l_{ij}=I_{ij}^2$; $I_{ij}$ is the current magnitude of line $ij$; and $\{\mathbb{L}\}$ is the set of transmission lines. Represented compactly, the ACOPF with *angle relaxation* is the following.

$$(\mathbf{AR}) \quad \min_{v,l,P^g,Q^g,P^c,Q^c,F,H} f(v,l,P^g,Q^g,P^c,Q^c,F,H) \quad (3\text{-a})$$

$$s.t. \quad (2\text{-a})\text{-}(2\text{-d}), \quad (v,\theta,l,P^g,Q^g,P^c,Q^c,F,H) \in \mathbb{X} \quad (3\text{-b})$$

where $\mathbb{X}$ represents individual variables' bound constraints if exist. Noting that quadratic equalities in (2-d) are non-convex, they can be relaxed as inequalities to have a conic relaxation.

$$(\mathbf{CR}) \quad \min_{v,l,P^g,Q^g,P^c,Q^c,F,H} f(v,l,P^g,Q^g,P^c,Q^c,F,H) \quad (4\text{-a})$$

$$s.t. \quad (2\text{-a})\text{-}(2\text{-c}), \quad (v,\theta,l,P^g,Q^g,P^c,Q^c,F,H) \in \mathbb{X} \quad (4\text{-b})$$

$$l_{ij}v_i \geq F_{ij}^2 + H_{ij}^2, \quad \forall ij \in \{\mathbb{L}\} \quad (4\text{-c})$$

It was shown in [1] that if the following four conditions or assumptions are met, **CR** derives an exact solution to **AR** and there is no gap between **AR** and **CR**.
(*i*) The network is connected and there are no upper bounds on $P^c$ and $Q^c$; (*ii*) The objective function $f$ is convex; (*iii*) $f$ is strictly increasing in $l$, non-increasing in $P^c$ and $Q^c$, and independent of $F$ and $H$; (*iv*) The OPF is a feasible model.

**Observations**: On an often used IEEE 300-bus system [5], we note that, although the aforementioned conditions are satisfied, conic constraints in (4-c) might not be tight and hence **CR** does not lead to an exact solution of **AR**. Hence, we believe that those conditions are not sufficient to guarantee the zero gap between **AR** and **CR**.

**Analysis**: To understand such discrepancy, we study the proof of exactness of **CR** presented in [1]. Consider an optimal solution ($v, l, P^g, Q^g, P^c, Q^c, F, H$) to **CR**. In [1], by contradiction, it is assumed that ($v, l, P^g, Q^g, P^c, Q^c, F, H$) is not optimal to AR. It is argued that for one branch $ij$ whose constraint in (4-c) is not binding, another solution constructed using (5) with a small ε > 0 is feasible (due to condition (*i*)) and has a better objective function value. Hence, the desired contradiction follows and ($v, l, P^g, Q^g, P^c, Q^c, F, H$) is optimal to **AR**.

Let $x'_{ij} = \frac{1}{2}x_{ij}\varepsilon$ and $r'_{ij} = \frac{1}{2}r_{ij}\varepsilon$.

$$\overline{v} = v, \overline{P}^g = P^g, \overline{Q}^g = Q^g, \overline{l}_{ij} = l_{ij} - \varepsilon, \overline{l}_{-ij} = l_{-ij} \quad (5\text{-a})$$

$$\overline{P}^c_i = P^c_i + r'_{ij}, \overline{P}^c_{-i} = P^c_{-i}, \overline{Q}^c_i = Q^c_i + x'_{ij}, \overline{Q}^c_{-i} = Q^c_{-i} \quad (5\text{-b})$$

$$\overline{P}^c_j = P^c_j + r'_{ij}, \overline{P}^c_{-j} = P^c_{-j}, \overline{Q}^c_j = Q^c_j + x'_{ij}, \overline{Q}^c_{-j} = Q^c_{-j} \quad (5\text{-c})$$

$$\overline{F}^c_{ij} = F_{ij} - r'_{ij}, \overline{F}^c_{-ij} = F_{-ij}, \overline{H}_{ij} = H_{ij} - x'_{ij}, \overline{H}_{-ij} = H_{-ij} \quad (5\text{-d})$$

Nevertheless, we notice that this argument [1] is actually based on an assumption that $x_{ij}$ is positive for this non-binding branch $ij$. If it has negative reactance, i.e., $x_{ij}<0$, the new solution might not be feasible, as according to (5-b) or (5-c), the new solution could violate the lower bound constraints on $\overline{Q}^c_i$ or $\overline{Q}^c_j$. Therefore, it could not lead to the desired contradiction.

Negative reactance may result from mutual effects between windings of the 3-winding transformers [3]. It may be caused by series compensation in transmission lines or Thevenin equivalent parameters for some parts of the power system. Actually, we note that branches with negative reactance exist in typical power systems. There is one such branch IEEE 300-bus system, 10 such branches in Polish 3012-bus and 3120 systems, and 12 such branches in Polish 3375-bus system [5]. Therefore, to ensure the correctness, the zero-gap sufficient conditions should be augmented to include an additional condition: *(v) there is no line with negative reactance*.

## II. NUMERICAL EXAMPLE

To evaluate the impact of negative reactance (NR), we first perform an experiment on IEEE 14-bus system by just changing one transmission line's reactance to its negative value, with the objective function of being the network losses. We then solve the **CR** model using CPLEX as the SOCP solver and compute the gap between the left-hand-side and the right-hand-side of conic constraint (4c) for each line. Results are listed in Table I where column "# of L/G" the number lines with non-zero gaps and column "Max Gap" represents the maximal gap. It can be seen from Table I that a single NR branch could seriously affect the relaxation quality of **CR**.

Next, we investigate IEEE 300-bus system. We compute ACOPF using interior point method with multiple starting points by IPOPT solver [4] as a benchmark. On line #178, we note a non-zero gap (0.2466) in the solution of the **CR** model, which results in a large discrepancy of voltage magnitudes at bus #248 to the benchmark method. It is because of both angle relaxation and conic relaxation. Figure 1 provides voltage magnitudes of all buses from these two ACOPF solutions. As mentioned, there is one transmission line with negative reactance, i.e., line #179 with the value of -0.3697p.u.. In contrast, we modify it to 0.3697p.u. and do the computation again. For new ACOPF solutions, we find that, as expected, all the conic constraints are binding and there is no gap between **AR** and **CR**. Fig. 2 shows voltage magnitudes of all buses, where there are marginal differences between the **CR** model (actually the AR model) and benchmark method. Obviously, it confirms the tightness of the angle relaxation **AR** model with respect to the original ACOPF formulation.

Table I   Impact of NE on the gap of conic constraints in model (4)

| Line # | Max gap | # of L/G | Line # | Max gap | # of L/G |
|---|---|---|---|---|---|
| 1 | 2.7329e-8 | -- | 11 | 1.9176e-8 | -- |
| 2 | 1.2494e-8 | -- | 12 | 5.6446e-7 | -- |
| 3 | 2.7329e-7 | -- | 13 | 2.3123e-8 | -- |
| 4 | 1.7781e-7 | -- | 14 | 28.1263 | 1 |
| 5 | 6.4705e-7 | -- | 15 | 1218.2 | 1 |
| 6 | 1.3354 | 1 | 16 | 1.7648 | 1 |
| 7 | 2.2503 | 1 | 17 | 0.8058 | 1 |
| 8 | 9.9505 | 1 | 18 | 1.4896 | 1 |
| 9 | 1.3994 | 1 | 19 | 7.0376e-8 | -- |
| 10 | 2.0293 | 1 | 20 | 2.4707e-7 | -- |

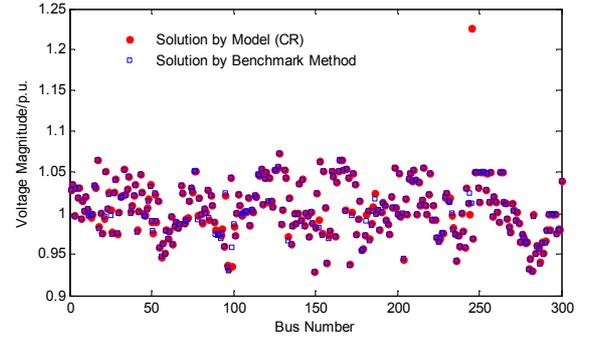

Fig.1 Comparison of **CR** with IPOPT benchmark with NR

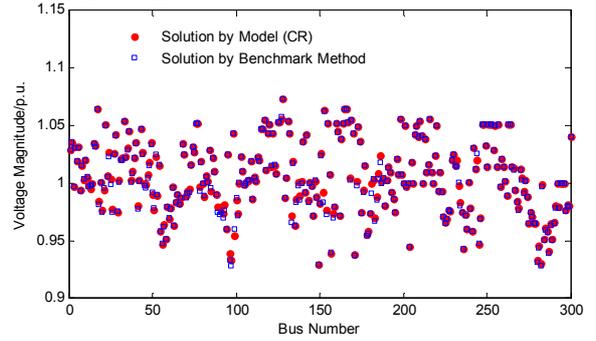

Fig.2 Comparison of **CR** with IPOPT benchmark without NR

## III. CONCLUSIONS

This letter shows that a practical issue, i.e., negative reactance, could cause the sufficient conditions presented in [1] fail to ensure the exactness of the conic relaxation. Hence, those sufficient conditions should be augmented to include that there is no line with negative reactance in the network. Numerical results confirm the necessity of this new condition.

REFERENCES

[1] M. Farivar, S. H. Low, "Branch Flow Model: Relaxations and Convexification-Part I," *IEEE Trans. Power Syst.*, vol. 28, no. 3, pp. 2554-2564, 2013.
[2] M. Farivar, S. H. Low, "M. Farivar, S. H. Low, "Branch Flow Model: Relaxations and Convexification-Part I,"," *IEEE Trans. Power Syst.*, vol. 28, no. 3, pp. 2565-2572, 2013.
[3] W. W. Weaver Jr, *Geometric and Game-theoretic Control of Energy Assets in Small-scale Power Systems*. ProQuest, 2007.
[4] A. Wächter and L. Biegler, "On the implementation of an interior-point filter line-search algorithm for large-scale nonlinear programming," *Math. Prog.*, vol. 106, pp. 25-57, 2006.
[5] R. D. Zimmerman, C. E. Murillo-Sánchez, and R. J. Thomas, "MATPOWER: Steady-State Operations, Planning, and Analysis Tools for Power Systems Research and Education," *IEEE Trans. Power Syst.*, vol. 26, no. 1, pp. 12-19, Feb 2011.